\documentclass{amsart}
\usepackage[T1]{fontenc}
\usepackage[utf8]{inputenc}
\usepackage{lmodern}
\usepackage{amssymb}
\usepackage[mathscr]{eucal}
\usepackage{tikz}
\usepackage{enumitem,booktabs}
\usetikzlibrary{arrows.meta,decorations.pathreplacing}
\newcommand{\deq}{\mathrel{\mathop:}=}
\newcommand{\alg}[1]{\mathbf{#1}}
\newcommand{\var}[1]{\mathsf{#1}}

\newcommand{\ov}[1]{\overline{#1}}
\newcommand{\Con}[1]{\mathbf{Con}\, #1}
\newcommand{\con}[1]{\mathrm{Con}\, #1}
\newcommand{\Cm}[1]{\mathbf{Cm}\, #1}
\newcommand{\cm}[1]{\mathrm{Cm}\, #1}

\newcommand{\Log}{\mathrm{Log}}
\def\ra{\rightarrow}
\def\Ra{\Rightarrow}
\def\up#1{{#1}\mathord{\uparrow}}
\def\dw#1{{#1}\mathord{\downarrow}}
\def\pw{\mathscr P}
\def\op{\bigstar}
\newcommand{\Ac}{\mathrm{Ac}}
\newcommand{\Nstar}{\mathbb N_{\geq 4}}
\newcommand{\doi}[1]{\href{https://doi.org/#1}{\nolinkurl{doi:#1}}}
\makeatletter
\newcommand{\doublewedge}{\@doubleop{\wedge}} 
\newcommand{\@doubleop}[1]{\DOTSB\mathop{\mathpalette\@doubleop@aux{#1}}}
\newcommand\@doubleop@aux[2]{%
  \sbox\z@{$\m@th#1#2$}%
  \makebox[1.25\wd\z@][s]{$\m@th#1#2\hss#2$}}
\makeatother
\newtheorem{theorem}{Theorem}[section]
\newtheorem{lemma}[theorem]{Lemma}
\newtheorem{cor}[theorem]{Corollary}
\newtheorem{prop}[theorem]{Proposition}
\theoremstyle{definition}
\newtheorem{defn}[theorem]{Definition}

\theoremstyle{remark}
\newtheorem{rem}[theorem]{Remark}
\tikzset{bdot/.style={circle,fill,inner sep=1.65pt},
  edge/.style={line width=.55pt},every label/.style={font=\small}}

\title[Rule extensions of the implication--negation fragment of IPC]{Rule extensions of the implication--negation fragment of intuitionistic logic}
\author[T. Kowalski]{Tomasz Kowalski}
\address{Department of Logic, Jagiellonian University, Krak\'ow, Poland}
\email{tomasz.s.kowalski@uj.edu.pl}
\author[K. S\l omczy\'nska]{Katarzyna S\l omczy\'nska}
\address{Department of Mathematics, University of the National Education Commission, Krak\'ow, Poland}
\email{irena.korwin-slomczynska@uken.krakow.pl}
\date{}

\begin{document}

\begin{abstract}
We present an algebraic study of extensions 
of the implication--negation fragment 
$\mathrm{IPC}_{\to,\neg}$ of intuitionistic propositional logic.
Such extensions by inference rules correspond to subquasivarieties of the
variety $\var{Ho}$ of Hilbert algebras with zero. We describe
the lattice of such extensions of $\mathrm{IPC}_{\to,\neg}$
by means of its generators: it is isomorphic to the lattice of upsets
of a countable poset, and in particular it is
completely distributive. Consequently,
$\var{Ho}$ is not Q-universal, in contrast to several other varieties
related to $\mathrm{IPC}$. Nevertheless, this
lattice is very complex: at each suitably defined height it
has intervals containing continuum many quasivarieties which are not varieties.
To the best of our knowledge $\mathrm{IPC}_{\to,\neg}$ is the only natural
fragment of intuitionistic logic known to combine a completely distributive
lattice of extensions with continuum many non-axiomatic extensions.
\end{abstract}
\maketitle

\section{Introduction}\label{sec:intro}

Hilbert algebras, introduced by Henkin~\cite{Hen50} under the name
\emph{implicative models}, and studied in depth by
Diego~\cite{Die66}, are the algebraic models of the implicational fragment
of intuitionistic propositional logic $\mathrm{IPC}$. Their expansions by a
least element $0$, which we call \emph{Hilbert algebras with zero}
(\emph{bounded Hilbert algebras} is another name in use), play the same
role for the fragment $\mathrm{IPC}_{\to,\neg}$ with implication and
negation, where $\neg x$ abbreviates $x\to0$. Both fragments are
algebraizable in the sense of Blok and Pigozzi~\cite{BP89} (see
also~\cite{Fon16}), with the variety $\var{Hi}$ of Hilbert algebras and the
variety $\var{Ho}$ of Hilbert algebras with zero as their equivalent
algebraic semantics. As a standard consequence of algebraizability,
the lattice of all finitary extensions
of $\mathrm{IPC}_{\to,\neg}$---by inference rules as well as by
axioms---is dually isomorphic to the lattice
$\Lambda^q(\var{Ho})$ of subquasivarieties of $\var{Ho}$, and under this
isomorphism the axiomatic extensions correspond to subvarieties. To be 
more precise, a quasivariety $\mathcal Q\subseteq\var{Ho}$
corresponds to the extension of $\mathrm{IPC}_{\to,\neg}$ by all rules
valid in every member of $\mathcal Q$, and $\mathcal{Q}$ is a variety if and
only if these rules are equivalent to zero-premiss rules, that is, to axioms.

The present article describes the structure $\Lambda^q(\var{Ho})$, pointing out 
sharp contrasts between it and subquasivariety lattices of other classes of
algebras related to $\mathrm{IPC}$. We work on the
algebraic side throughout; the dual isomorphism above gives
a translation of our results to logics, which we spell out
in Sections~\ref{subsec:logic1} and~\ref{subsec:logic2}.

Prucnal~\cite{Pru72} shows that the implicational fragment $\mathrm{IPC}_{\to}$ of
$\mathrm{IPC}$ is structurally complete: every admissible rule is
derivable. The same notion (in a suitable equivalent formulation,
see~\cite{Ber91}) is used in algebra, so $\var{Hi}$ is structurally complete, 
in fact \emph{hereditarily structurally complete} in the sense
of Citkin~\cite{Cit78}, that is, every
subquasivariety of $\var{Hi}$ is a variety. Minari and Wro\'nski~\cite{MW88}
showed that the same holds for
$\mathrm{IPC}_{\to,\wedge}$ and $\mathrm{IPC}_{\to,\wedge,\neg}$.
The implication-negation fragment behaves differently:
Wro\'nski~\cite{Wro86} showed that $\mathrm{IPC}_{\to,\neg}$ is not
structurally complete, Cintula and Metcalfe~\cite{CM10} gave a basis for
its admissible rules, and unification and projectivity in this setting were
studied in~\cite{Slo12}; for admissibility in full $\mathrm{IPC}$ we refer
to Rybakov~\cite{Ryb97} and Iemhoff~\cite{Iem01}. In algebraic terms, lack of
structural completeness means that 
$\var{Ho}$ has subquasivarieties which are not varieties: out of
the four varieties studied in~\cite{KSfree}---Hilbert algebras and Brouwerian
semilattices, as well as their expansions by zero---it is the only one with this
property. It is therefore natural to investigate
$\Lambda^q(\var{Ho})$ closer.

Our first main result, in Section~\ref{sec:distributive}, is the description of
the completely join-prime elements of $\Lambda^q(\var{Ho})$, which yields
a reasonably complete
description of the whole $\Lambda^q(\var{Ho})$: it is isomorphic to the
lattice of downsets of the poset of its completely join-prime elements,
a countable poset with a decidable order. In particular, $\Lambda^q(\var{Ho})$
is completely distributive, which 
rules out \emph{Q-universality} of $\var{Ho}$. The notion was
introduced by Sapir~\cite{Sap85}: a quasivariety $\mathcal{K}$ of finite
signature is Q-universal if the
subquasivariety lattice of every quasivariety of finite signature is a
homomorphic image of a sublattice of $\Lambda^q(\mathcal K)$; such lattices
satisfy no nontrivial lattice identity (see Adams and Dziobiak~\cite{AD94}).
The variety of Heyting algebras $\var{H}$, the variety of
$\{\wedge,\vee,\neg,0,1\}$-subreducts of $\var{H}$ (with $\neg x\deq x\ra 0$),
known as \emph{distributive p-algebras} or simply \emph{p-algebras},
and the variety of $\{\wedge,\neg\}$-subreducts of $\var{H}$,
known as \emph{pseudocomplemented semilattices},
are Q-universal (see~\cite{Dzi85}, \cite{AD94} and~\cite{ADGS95}).
Known classes of
subreducts of Heyting 
algebras that are not Q-universal include Brouwerian semilattices, Brouwerian
semilattices with zero, and Hilbert algebras: these are not Q-universal since
all their subquasivarieties are varieties, and so their sub(quasi)variety lattices
are distributive (which in turn follows from their congruence distributivity and
J\'onsson Lemma). 

The second main result, in Section~\ref{sec:continua}, is
that intervals of continuum size of proper subquasivarieties
occur at each finite \emph{height} (see~\ref{subsec:heights} below for the
definition). As far as we were able to determine,
$\var{Ho}$ is the first identified class of subreducts of Heyting algebras
which is (i) not Q-universal and (ii) has a continuum of subquasivarieties that
are not varieties. We note that the implication--disjunction fragment
$\mathrm{IPC}_{\to,\vee}$ is not structurally complete either, as
Mints' rule uses only $\to$ and $\vee$, but we are not aware of any
description of its lattice of extensions by rules.

\section{Preliminaries}\label{sec:prelim}

We rely on the results and notation of~\cite{KSfree},
in particular on the antichain representation of finite Hilbert algebras.
For completeness we restate some basic facts and calculations
needed for the present construction, but it
would be convenient for the reader to consult~\cite{KSfree} frequently.
We assume some knowledge of universal algebra and lattice theory:
a few chapters of~\cite{Ber12} suffice for the former, a few chapters
of~\cite{Gra11} for the latter.

\subsection{Hilbert algebras and congruences}
An algebra $\alg A=(A;\ra,1)$ is a Hilbert algebra if it satisfies
\begin{align*}
  x\ra x&=1, \\
  1\ra x&=x,\\
 x\ra(y\ra z)&=(x\ra y)\ra(x\ra z),\\
 (x\ra y)\ra((y\ra x)\ra x)&=(y\ra x)\ra((x\ra y)\ra y).
\end{align*}
Smaller equational bases exist, but this one is reasonably natural, and moreover
has been used in~\cite{KSfree}. Every Hilbert algebra is ordered by the relation:
$x\leq y$ if and only if $x\ra y=1$. The variety $\var{Ho}$ consists of
the expansions $(A;\ra,1,0)$ satisfying $0\ra x=1$.
The two identities below hold in every Hilbert algebra
\begin{equation}\label{eq:exchange}
 x\ra(y\ra z)=y\ra(x\ra z),\qquad
 x\ra(x\ra z)=x\ra z
\end{equation}
and justify the following shorthand notation.
For a finite set of terms $\{t_1,\ldots,t_k\}$, we write
$\{t_1,\ldots,t_k\}\ra r$ for 
$t_1\ra(t_2\ra\cdots\ra(t_k\ra r)\cdots)$, as in~\cite{KSfree}.
By~\eqref{eq:exchange} we can permute antecedents, so 
the order of the $t_i$ is irrelevant; hence the set notation.

We write $\Con{\alg A}$ for the congruence lattice and $\Cm{\alg A}$
for its poset of completely meet-irreducible congruences, ordered by
inclusion. The congruence generated by $(a,b)$ is denoted
$\theta^{\alg A}(a,b)$, with the superscript omitted when unambiguous.
Hilbert algebras are locally finite\footnote{An argument proving local
  finiteness of Hilbert algebras, whose variants were used later in many
  local finiteness proofs, appeared in Diego~\cite{Die66}.},
congruence distributive
and congruence $1$-regular: a congruence is determined
by its $1$-class. More explicitly, for every $\theta\in\con{\alg A}$,
\begin{equation}\label{eq:one-regular}
 (u,v)\in\theta\ \Longleftrightarrow\
 (u\ra v,1)\in\theta\text{ and }(v\ra u,1)\in\theta,
\end{equation}
and in particular for principal congruences the following characterisation holds:
\begin{equation}\label{eq:principal}
 (u,v)\in\theta(1,a)\ \Longleftrightarrow\ a\ra u=a\ra v.
\end{equation}

It is well known that a Hilbert algebra is subdirectly irreducible
if and only if it has a greatest element strictly below $1$.
We denote this element by $\op$, as in~\cite{KSfree}.

\subsection{Algebras of antichains}
We recall the construction from \cite[Section~5]{KSfree}.
For a finite poset $P=(P;\leq)$, let $\Ac(P)$ be the set of all
antichains of $P$. For $C\subseteq P$ put
$\dw C=\{p\in P:p\leq c\text{ for some }c\in C\}$.
For $B,C\in\Ac(P)$ define
\begin{equation}\label{eq:antichains}
 C\preccurlyeq B\ \Longleftrightarrow\ B\subseteq\dw C,
 \qquad
 C\ra B\deq B\setminus\dw C
       =\{b\in B:\up b\cap C=\emptyset\}.
\end{equation}
The order on the singletons is the reverse of that on $P$.
The algebra $\alg{Ac}(P)=(\Ac(P);\ra,\emptyset)$ belongs to
$\var{Hi}$, and its least element is $M\deq\max P$.

For comparison with the usual construction of Heyting algebras (and their
various reducts) as algebras of upsets, consider the bijection
\[
 \varphi:\Ac(P)\longrightarrow\mathrm{Up}(P),\qquad
 \varphi(C)=P\setminus\dw C.
\]
Now, implication on upsets is given by
$U\Rightarrow V=\{p:\up p\cap U\subseteq V\}$, so we have
\[
 \varphi(C)\Rightarrow\varphi(B)
 =P\setminus\dw{(\dw B\setminus\dw C)}
 =P\setminus\dw{(B\setminus\dw C)}
 =\varphi(C\ra B).
\]
Hence, $\alg{Ac}(P)$ is isomorphically mapped to the implication
reduct of the Heyting algebra of upsets of $P$; in particular,
$\varphi(M) = \emptyset$ is the abstract element $0$,
and $\varphi(\emptyset) = P$ is the abstract element $1$. We will not use the
upset algebras any more, but we will keep using the abstract $0$ and $1$
whenever convenient, so the reader should not be alarmed by an occasional
occurrence of $\emptyset = 1$ or similar.

To obtain Hilbert algebras we only need to take certain
  subfamilies of $\Ac(P)$. Let us call a family $S\subseteq \Ac(P)$
  \emph{hereditary} if $S$ is closed under taking subsets. Note that
a nonempty hereditary $S$ always has $\emptyset\in S$.

\begin{lemma}\label{lem:selection}
Let a nonempty $S\subseteq\Ac(P)$ be hereditary.
Then $(S;\ra,\emptyset)$ is a subalgebra of
$\alg{Ac}(P)$. If $\max P\in S$, the expansion
$(S;\ra,\emptyset,\max P)$ belongs to $\var{Ho}$.
\end{lemma}
\begin{proof}
The value $C\ra B$ is a subset of $B$. Hence it belongs to $S$ whenever
$B,C\in S$. If $M=\max P$ belongs to $S$, finiteness gives
$\dw M=P$, so $M\ra B=\emptyset$ for every $B\in S$.
\end{proof}

Conversely, every finite member of $\var{Ho}$ arises in this way
as~\cite[Theorems~5.6--5.8]{KSfree}
shows. We recall the result below for easier reference.  

\begin{theorem}\label{thm:antichains}
Each finite algebra $\alg{A}\in\var{Ho}$ is isomorphic to a 
subalgebra of the zero-expansion
$(\Ac(P);\ra,\emptyset,\max P)$ of $\alg{Ac}(P)$,
whose universe is a hereditary
subset of $\Ac(P)$ containing all singletons,
where $P = (\cm{\alg{A}};\subseteq)$.  
\end{theorem}

\begin{defn}\label{def:represented-over}
Let $P$ be a finite poset. A finite algebra $\alg{A}\in\var{Ho}$ is
\emph{represented over} $P$ if
$\alg{A}\cong (S;\ra,\emptyset,\max P)$ for some
hereditary $S\subseteq \Ac(P)$ with $\max P\in S$ and
$\{p\}\in S$ for all $p\in P$.
\end{defn}  

By Theorem~\ref{thm:antichains}, every finite member of $\var{Ho}$ is
represented over some finite poset. The next lemma shows that this poset is
determined by the algebra up to isomorphism.

\begin{lemma}\label{lem:height-poset}
Let $\alg A$ be represented over $P$. Then
$\Cm{\alg A}\cong P$.
\end{lemma}

\begin{proof}
We may assume $\alg A=(A;\ra,\emptyset,\max P)$ with $A\subseteq\Ac(P)$
as in Definition~\ref{def:represented-over}.
Since $\Con{\alg A}$ is a finite distributive lattice, its posets of
join-irreducible and meet-irreducible elements are isomorphic
(see, e.g., \cite[Chapter~II, Section~1]{Gra11}).
We write $J(\Con{\alg A})$ for the
poset of join-irreducible elements of $\Con{\alg A}$;
by \cite[Theorem~5.1]{KSfree}, its elements are exactly
the congruences $\theta(1,U)$ with $U\neq\emptyset$ such that
\begin{equation}\label{eq:irreducible-element}
 C\ra U\in\{\emptyset,U\}\qquad\text{for every }C\in A.
\end{equation}
Every singleton $U=\{p\}$ has this property, but if $|U|\geq2$, then for
$p\in U$ we have $\{p\}\ra U=U\setminus\{p\}$, which is neither
$\emptyset$ nor $U$. Consequently,
\[
 J(\Con{\alg A})=\{\theta(1,\{p\}):p\in P\}.
\]
Moreover, by~\eqref{eq:principal},
\[
 \theta(1,\{p\})\subseteq\theta(1,\{q\})
 \quad\Longleftrightarrow\quad
 \{q\}\ra\{p\}=\emptyset
 \quad\Longleftrightarrow\quad p\leq q.
\]
Thus $p\mapsto\theta(1,\{p\})$ is an order isomorphism from $P$
onto $J(\Con{\alg A})$. Therefore, $\Cm{\alg A}\cong P$.
\end{proof}

We shall also use the following concrete description of quotients.
Readers familiar with Kripke frames will recognise generated subframes there.
Indeed, the representation of $\alg{A}$ over $P$ can be seen as a Kripke frame
with a selected set of antichains. 

\begin{lemma}\label{lem:restriction}
Let $U$ be an upset of a finite poset $P$, let $M=\max P$, and let
$\alg S=(S;\ra,\emptyset,M)$ be a subalgebra of the zero-expansion of
$\alg{Ac}(P)$. Put $S_U\deq\{C\cap U:C\in S\}$ and let
$\alg{S}_U\deq (S_U;\ra,\emptyset,M\cap U)$, with implication computed in $U$.
Then $r_U\colon S\to S_U$ given by
\[
r_U(C)\deq C\cap U
\]
is a surjective homomorphism from $\alg S$ onto $\alg{S}_U$ and its
kernel is the congruence
$\{(B,C)\in S^2:B\cap U=C\cap U\}$.
Moreover, if $S$ is hereditary and $M\subseteq U$, then $\alg S_U$ is a
subalgebra of $\alg S$ and $r_U$ is an idempotent endomorphism of $\alg S$;
thus $\alg S_U$ is a retract of $\alg S$.
\end{lemma}

\begin{proof}
For $b\in U$, if $b\leq c$ in $P$, then $c\in U$ because $U$ is an
upset. Thus
\[
r_U(C\ra B) =  (B\setminus\dw C)\cap U
=(B\cap U)\setminus\bigl((C\cap U)\mathord{\downarrow}_{U}\bigr)
= (C\cap U) \ra (B\cap U),
\]
where $\mathord{\downarrow}_{U}$ denotes the downset in $U$. This proves
preservation
of implication. Both constants are preserved since
$M\cap U=\max U$. The description of the kernel is immediate.

To show the moreover part, note that $r_U$ is idempotent by definition,
and by $M\subseteq U$ we have $S_U\subseteq S$. It is easy to see that
$S_U$ is closed under implication, so
$\alg{S}_U$ is a subalgebra of $\alg{S}$. Hence,
$r_U\colon S\to S$ is a retraction, and so
$\alg S_U=r_U(\alg S)$ is a retract of $\alg S$.
\end{proof}

\begin{rem}\label{rem:boolean}
Let $S$ be hereditary and take $U=M$ in Lemma~\ref{lem:restriction}.
Since the points of $M$ are maximal, for every $C\in S$ we have
$\neg C=C\ra M=M\setminus\dw C=M\setminus C$, and therefore
$ r_M(C)=C\cap M=M\setminus(M\setminus C)=\neg\neg C$,
that is, $r_M$ is the double negation map, which is known to be
a retraction for any Hilbert algebra with zero.
Since $S$ is hereditary and contains $M$, the image $r_M(S)$ is  
$(\pw(M); \ra, \emptyset, M)$, clearly isomorphic to
$\alg 2^{|M|}$ (cf.~Section~\ref{subsec:heights}). Note also that
$\alg 2$ is isomorphic to 
the subalgebra $\{0,1\}$ of any nontrivial member of $\var{Ho}$; hence
$\alg 2$, together with all its powers, belongs to every nontrivial
subquasivariety of $\var{Ho}$. These facts will be used in
Section~\ref{sec:distributive}.
\end{rem}

\subsection{Quasivarieties}\label{subsec:quasivarieties}

We shall use a few facts of quasivariety-theory folklore. All background on
quasivarieties and subquasivariety lattices we need here (and much more) can
be found in Hyndman and Nation~\cite{HN18}. As usual,
$I,S,P,P_{\mathrm{fin}}$ and $P_U$ stand for class operators
of taking isomorphic copies, subalgebras, direct products, finite direct products,
and ultraproducts, respectively. Finite algebras that do not embed into products of
their proper subalgebras are usually called \emph{critical}, but 
they are called \emph{quasicritical} in~\cite{HN18}. The next lemma shows
that, conversely, an algebra which does embed into a product of its proper
subalgebras is dispensable as a generator of a quasivariety.

\begin{lemma}\label{lem:generating}
Let $\alg{A}$ be a subalgebra of a product $\prod_{i\in I}\alg{B}_i$ such that
$\alg{B}_i$ is a subalgebra of $\alg{A}$ for every $i\in I$. Then
$Q(\alg{A}) = Q(\{\alg{B}_i: i\in I\}) = \bigvee_{i\in I} Q(\alg{B}_i)$.
\end{lemma}

\begin{lemma}\label{lem:finite-products}
For a finite algebra $\alg A$ in a finite signature, and any class
$\mathcal K$ of algebras of that signature,
\[
 \alg A\in Q(\mathcal K)
 \quad\Longleftrightarrow\quad
 \alg A\in ISP_{\mathrm{fin}}(\mathcal K).
\]
\end{lemma}
The reader interested in a proof may consult~{\cite[Lemma~1.6]{KSpa}}.

\begin{cor}\label{cor:one-factor}
If $\alg A$ is finite and subdirectly irreducible, in a finite signature,
then
\[
 \alg A\in Q(\mathcal K)\quad\Longleftrightarrow\quad
 \alg A\in IS(\mathcal K).
\]
\end{cor}

\begin{proof}
Since $\alg{A}\in ISP_{\mathrm{fin}}(\mathcal K)$ is subdirectly irreducible,
it must embed into one of its factors, hence
$\alg{A}\in IS(\mathcal K)$. The converse is obvious. 
\end{proof}

\section{The lattice of subquasivarieties}\label{sec:distributive}

Throughout this section, $\alg A$ will be an algebra represented
over a finite poset $P$, with $M=\max P$. An essential role will be played
by algebras that are, intuitively, not far from being subdirectly irreducible.

\begin{defn}\label{def:nearly-si}
A finite algebra $\alg{A}\in\var{Ho}$ represented over a nonempty $P$ will be called
\emph{nearly subdirectly irreducible} if either $P = M$ or
$P\setminus M$ has a least element.
\end{defn}  

The concept of being nearly subdirectly irreducible could easily be extended to
apply also to infinite algebras, however, since we will never apply it so in
the present article we decided to include finiteness in the definition.
Clearly, every finite subdirectly irreducible member of $\var{Ho}$ is
nearly subdirectly irreducible. The trivial algebra, which is
represented over the empty poset, is not nearly subdirectly irreducible.

\begin{lemma}\label{lem:embedding}
Let $\alg A$ represented over $P$ be nearly subdirectly irreducible
with $P\neq M$ and $d \deq\min(P\setminus M)$.
For every homomorphism $f\colon\alg A\to\alg B$ into some
$\alg B\in\var{Ho}$ such that $f(\{d\})\neq 1$,
the map $\ov f\colon A\to B\times\pw(M)$ given by
\[
 \ov f(C)\deq\bigl(f(C),\,C\cap M\bigr)
\]
is an embedding of $\alg A$ into $\alg B\times\alg{2}^{|M|}$.
Hence, $\alg{A}\in Q(\alg{B})$. 
\end{lemma}

\begin{proof}
Note first that by Remark~\ref{rem:boolean}, the map
$C\mapsto C\cap M$ is a homomorphism from $\alg A$ onto
$(\pw(M);\ra,\emptyset,M)\cong \alg{2}^{|M|}$.
Thus both coordinates of $\ov f$ are homomorphisms, so $\ov f$ is. Since Hilbert
algebras are congruence $1$-regular, $\ov f$ is injective if and only if
$\ov f(C)=1$ implies $C=1=\emptyset$. Suppose then that $C\in A$ is
nonempty and $\ov f(C)=(1,\emptyset)$. Then $C\cap M=\emptyset$, so every
$c\in C$ lies in $P\setminus M$, and hence $d\leq c$. Thus
$\{d\}\subseteq\dw C$, which by~\eqref{eq:antichains} means
$C\preccurlyeq\{d\}$, or equivalently, $C\ra\{d\}=\emptyset=1$. Applying
$f$ we obtain
\[ 
 1=f(C)\ra f(\{d\})=1\ra f(\{d\})=f(\{d\}),
\]
contrary to the assumption. Hence, $\ov f$ is injective, and 
we have $\alg A\in IS(\alg B\times\alg 2^{|M|})$. 
As $\alg{B}$ is nontrivial, 
$\alg{B}\times\alg{2}^{|M|}\in Q(\alg{B})$ 
and therefore $\alg{A}\in Q(\alg{B})$.
\end{proof}

\begin{lemma}\label{lem:join-prime}
Let $\alg A$ represented over $P$ be nearly subdirectly irreducible
with $P\neq M$. For any family
$\{\mathcal Q_i:i\in I\}$ of subquasivarieties of $\var{Ho}$, the following holds.
If $\alg A\in\bigvee_{i\in I}\mathcal Q_i$, then $\alg A\in\mathcal Q_i$
for some $i\in I$. Hence, $Q(\alg A)$ is a completely join-prime
element of $\Lambda^q(\var{Ho})$.
\end{lemma}

\begin{proof}
Since $\alg A$ is finite and
$\bigvee_{i\in I}\mathcal Q_i=Q\bigl(\bigcup_{i\in I}\mathcal Q_i\bigr)$,
Lemma~\ref{lem:finite-products} gives an embedding
$f\colon\alg A\hookrightarrow\alg K_1\times\cdots\times\alg K_n$ with
$\alg K_j\in\mathcal Q_{i_j}$ for some $i_1,\ldots,i_n\in I$.
Let $d = \min(P\setminus M)$. There is a $j$ such that $\pi_j(f(\{d\}))\neq 1$,
where $\pi_j$ denotes the $j$-th projection; in particular, $\alg K_j$ is
nontrivial. Applying Lemma~\ref{lem:embedding} to the homomorphism
$\pi_j\circ f\colon\alg A\to\alg K_j$ we get 
$\alg A\in IS(\alg K_j\times\alg 2^{|M|})$ and
$\alg{A}\in Q(\alg{K}_j)$. Hence,
$\alg A\in\mathcal Q_{i_j}$. Therefore,
$Q(\alg{A})\subseteq \bigvee_{i\in I}\mathcal Q_i$ implies
$Q(\alg{A})\subseteq \mathcal Q_{i_j}$, and so $Q(\alg A)$ is completely
join-prime.
\end{proof}

To state the next lemma, some notation will be convenient. 
Let $\alg A$ be finite, and represented over $P$ with $M=\max P$.
Let $D$ be the set of minimal elements of
$P\setminus M$. If $D\neq\emptyset$, for each $d\in D$ put
\[
U_d\deq\up d\cup M,\qquad
 A_d\deq\{C\cap U_d: C\in A\},\qquad
 \alg A_d\deq(A_d;\ra,\emptyset,M),
\] 
so that $\alg A_d=\alg A_{U_d}$ in the notation of
Lemma~\ref{lem:restriction}. Note that $A_d$ is hereditary and contains
$M=\max U_d$ and all singletons $\{p\}$ with $p\in U_d$, so $\alg A_d$
is represented over $U_d$; moreover $d$ is the least element of
$U_d\setminus M$, so $\alg A_d$ is nearly subdirectly irreducible.

\begin{lemma}\label{lem:decomposition}
Let $\alg A$ be represented over $P$, with $M=\max P$.
\begin{enumerate}
\item If $P=M$, then $\alg A\cong\alg 2^{|M|}$.
\item If $P\neq M$, then each $\alg A_d$ is a retract of $\alg A$ and
$\alg{A}$ embeds into $\prod_{d\in D}\alg A_d$.
\end{enumerate}
\end{lemma}

\begin{proof}
(1) If $P=M$, then $\Ac(P)=\pw(M)$, so $A=\pw(M)$, because $A$ is
hereditary and contains $M$. Hence $\alg A\cong\alg 2^{|M|}$ by
Remark~\ref{rem:boolean}.

(2) Let $P\neq M$, so that $D\neq\emptyset$ by finiteness, and fix $d\in D$.
The set $U_d$ is an upset of $P$ containing $M$.
Hence, by Lemma~\ref{lem:restriction},
the map $r_d\colon A\to A$ given by $r_d(C)=C\cap U_d$ is an idempotent
endomorphism of $\alg A$, and so $\alg A_d$ is a retract of $\alg A$.
The product map
$\prod_{d\in D} r_d: \alg{A}\to \prod_{d\in D}\alg{A}_d$ is clearly a
homomorphism, and injectivity follows from the fact that $\bigcup_{d\in D}U_d = P$.
\end{proof}

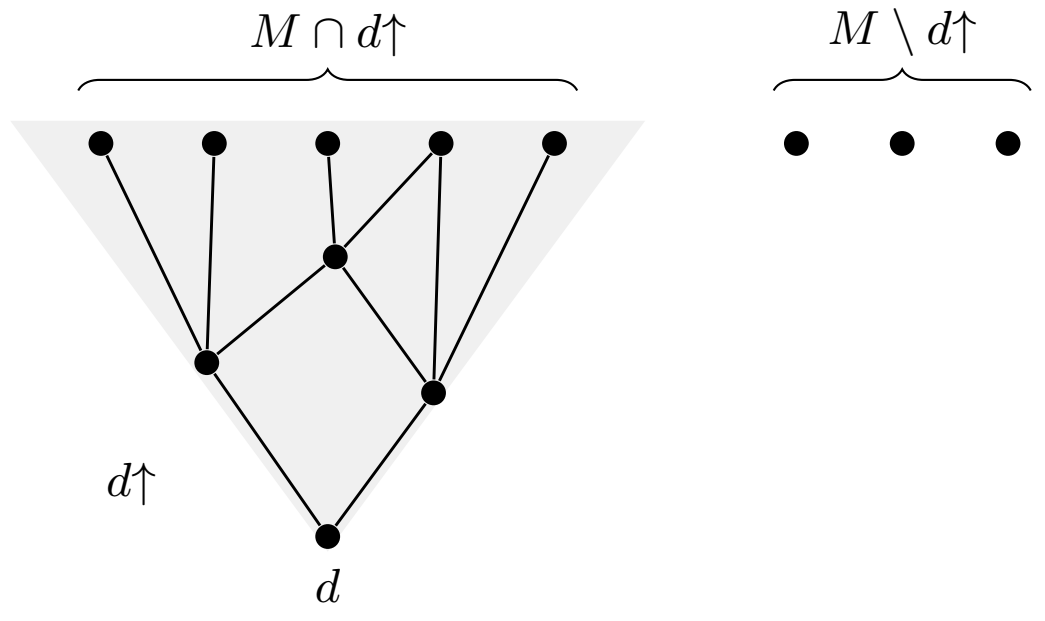
\begin{figure}[htbp]
\centering
\begin{tikzpicture}[x=1cm,y=1cm]
 \fill[gray!12] (0,-.1) -- (-2.1,2.75) -- (2.1,2.75) -- cycle;
 \node[bdot,label=below:{$d$}] (d) at (0,0) {};
 \node[bdot] (x1) at (-0.8,1.15) {};
 \node[bdot] (x2) at (0.7,0.95) {};
 \node[bdot] (x3) at (0.05,1.85) {};
 \foreach \x/\n in {-1.5/a,-0.75/b,0/c,0.75/e,1.5/f} \node[bdot] (m\n) at (\x,2.6) {};
 \draw[edge] (d)--(x1); \draw[edge] (d)--(x2);
 \draw[edge] (x1)--(x3); \draw[edge] (x2)--(x3);
 \draw[edge] (x1)--(ma); \draw[edge] (x1)--(mb);
 \draw[edge] (x3)--(mc); \draw[edge] (x3)--(me);
 \draw[edge] (x2)--(me); \draw[edge] (x2)--(mf);
 \node[font=\small] at (-1.3,0.35) {$\up d$};
 \draw[decorate,decoration={brace,amplitude=4pt}] (-1.65,2.95) -- (1.65,2.95)
   node[midway,above=3pt,font=\small] {$M\cap\up d$};
 \foreach \x in {3.1,3.8,4.5} \node[bdot] at (\x,2.6) {};
 \draw[decorate,decoration={brace,amplitude=4pt}] (2.95,2.95) -- (4.65,2.95)
   node[midway,above=3pt,font=\small] {$M\setminus\up d$};
\end{tikzpicture}
\caption{The poset $U_d=\up d\cup M$ representing the retract $\alg A_d$.}
\label{fig:Ad}
\end{figure}

\begin{rem}\label{rem:Ad}
Since $U_d$ of Lemma~\ref{lem:decomposition} is the disjoint union of
$\up d$ and the set $M'\deq M\setminus\up d$, the algebra $\alg{A}_d$ is a
subdirect product of 
$\alg{A}_{\up d}\deq(A_{\up d};\ra,\emptyset,M\cap\up d)$ 
where $A_{\up d}\deq\{C\in A:C\subseteq\up d\}$ and
$(\pw(M');\ra,\emptyset,M')\cong\alg 2^{|M'|}$.
Clearly, $\alg{A}_{\up d}$ is subdirectly irreducible,
with $\op = \{d\}$. We warn readers familiar with Kripke frames that despite
$U_d = \up{d}\uplus M'$ being a disjoint union, $\alg{A}_d$ is in general
not a direct product $\alg{A}_{\up d}\times \alg 2^{|M'|}$. Examples
illustrating this will play a crucial role in Section~\ref{sec:continua}.
\end{rem}

\begin{lemma}\label{lem:prime-singly-gen}
Let $\mathcal{Q}\in \Lambda^q(\var{Ho})$ be completely join-prime. Then
$\mathcal{Q} = Q(\alg{A})$ for a nearly subdirectly irreducible algebra
$\alg{A}\in\var{Ho}$. 
\end{lemma}  

\begin{proof}
Take a completely join-prime $\mathcal{Q}$. If $\mathcal{Q}$ is the smallest
nontrivial member of $\Lambda^q(\var{Ho})$, term-equivalent to the variety of
Boolean algebras, then $\mathcal{Q} = Q(\alg{2})$ and the claim holds.

Assume $\mathcal{Q}$ strictly contains $Q(\alg{2})$, so that
$\mathcal{Q} = Q(\{\alg{A}_i: i\in I\})$ for some family of   
non-Boolean $\alg{A}_i\in \mathcal{Q}$. Clearly,
$\mathcal{Q} = \bigvee_{i\in I} Q(\alg{A}_i)$, so by complete join-primeness
$\mathcal{Q} = Q(\alg{A}_i)$ for some $i\in I$. Without loss of generality,
assume $\mathcal{Q} = Q(\alg{A})$ for some $\alg{A}\in\var{Ho}$.
Using the facts that every algebra embeds into an ultraproduct of its finitely
generated subalgebras, and that $\var{Ho}$ is locally finite,
we quickly get that $Q(\alg{A})$ is generated by a set $\{\alg{B}_j:j\in J\}$ of
finite algebras, and so $Q(\alg{A}) = \bigvee_{j\in J}Q(\alg{B}_j)$.  
Now, since $Q(\alg{A})$ is completely join-prime,
$Q(\alg{A}) = Q(\alg{B}_j)$ for some $j\in J$.
Hence, $\mathcal{Q} = Q(\alg{B})$ for a finite $\alg{B}\in\var{Ho}$.

As $Q(\alg{B})=\mathcal{Q}\supsetneq Q(\alg{2})$, the algebra
$\alg{B}$ is not a Boolean algebra, so by 
Lemma~\ref{lem:decomposition}(2) it embeds into
a product $\prod_{d\in D}\alg{B}_d$ with each $\alg{B}_d$ a nearly
subdirectly irreducible retract of $\alg{B}$. Hence, by
Lemma~\ref{lem:generating} we get $Q(\alg{B}) = \bigvee_{d\in D} Q(\alg{B}_d)$
and join-primeness of $Q(\alg{B})$ finally yields that
$\mathcal{Q} = Q(\alg{B}) = Q(\alg{B}_d)$ for some  $d\in D$,
proving the claim.
\end{proof}

Lemmas~\ref{lem:join-prime} and~\ref{lem:prime-singly-gen} immediately yield the
result below.

\begin{cor}\label{cor:join-prime}
For any quasivariety $\mathcal{Q}\in\Lambda^q(\var{Ho})$ the following ere
equivalent:
\begin{enumerate}
\item $\mathcal{Q}$ is completely join-prime.
\item $\mathcal{Q} = Q(\alg{A})$ for a nearly subdirectly irreducible
  $\alg{A}\in\var{Ho}$.  
\end{enumerate}  
\end{cor}  

Now we are ready to prove our first theorem on the structure of
$\Lambda^q(\var{Ho})$.  

\begin{theorem}\label{thm:structure}
Let\/ $\mathbb{P} = \big(\{Q(\alg{A})\in \Lambda^q(\var{Ho}): \alg{A} \text{ nearly
  subdirectly irreducible}\}; \subseteq\big)$.  Then
$$
\Lambda^q(\var{Ho}) \cong \mathcal{O}(\mathbb{P})
$$
where $\mathcal{O}(\mathbb{P})$ is the lattice of downsets of\/ $\mathbb{P}$.
\end{theorem}  

\begin{proof}
For any $\mathcal{Q}\in\Lambda^q(\var{Ho})$ let
$\mathbb{N}_{\mathcal{Q}}$ be the set of all nearly subdirectly irreducible
$\alg{A}\in\mathcal{Q}$; since nearly subdirectly irreducible algebras are
nontrivial, $\mathbb{N}_{\mathcal{Q}}=\emptyset$ when $\mathcal{Q}$ is the
trivial quasivariety.
By Lemmas~\ref{lem:generating} and~\ref{lem:decomposition}(2) we have 
$\mathcal{Q} = Q(\mathbb{N}_{\mathcal{Q}})$.
Define a map $\Phi$ from $\Lambda^q(\var{Ho})$ to $\mathcal{O}(\mathbb{P})$
putting
$$
\Phi(\mathcal{Q})\deq
\begin{cases}
\emptyset & \text{ if } \mathcal{Q} \text{ is trivial},\\ 
\{Q(\alg{A})\subseteq\mathcal{Q}:\alg{A}\in\mathbb{N}_{\mathcal{Q}}\} &
\text{ otherwise}.
\end{cases}                                                                        
$$
This is clearly a downset of $\mathbb{P}$, and $\Phi$ is order-preserving.
It is also order-reflecting, hence injective: if
$\Phi(\mathcal{Q}_1)\subseteq\Phi(\mathcal{Q}_2)$, then
$\mathbb{N}_{\mathcal{Q}_1}\subseteq\mathcal{Q}_2$, and so
$\mathcal{Q}_1=Q(\mathbb{N}_{\mathcal{Q}_1})\subseteq\mathcal{Q}_2$.
For surjectivity, let $D$ be a downset of $\mathbb{P}$ and
let $\mathcal{Q} = \bigvee D$. Take any $\alg{A}\in \mathbb{N}_{\mathcal{Q}}$; then
$\alg{A}\in\mathcal{R}$ for some $\mathcal{R}\in D$ by
Lemma~\ref{lem:join-prime}. Hence,  
$Q(\alg{A})\subseteq\mathcal{R}$, and as
$D$ is a downset, $Q(\alg{A})\in D$, so $\Phi(\mathcal{Q})\subseteq D$. The
converse inclusion is immediate by definition of $D$, so $\Phi(\mathcal{Q}) =
D$ as required.
\end{proof}

\begin{cor}\label{cor:distributive}
The lattice $\Lambda^q(\var{Ho})$ is completely distributive.
Consequently, $\var{Ho}$ is not $Q$-universal.
\end{cor}

\begin{proof}
Let $\mathbb{P}$ be as in Theorem~\ref{thm:structure}. 
Since the lattice operations 
on $\mathcal{O}(\mathbb{P})$ are set-theoretical,
complete distributivity of $\Lambda^q(\var{Ho})$ follows. In particular,
$\Lambda^q(\var{Ho})$ is not Q-universal.  
\end{proof}

Since $\mathbb{P}$ determines the structure of
$\Lambda^q(\var{Ho})$, it is useful to observe that
the inclusion order on $\mathbb{P}$ is decidable.

\begin{theorem}\label{thm:inclusion-charact}
Let $\alg{A}$ and $\alg{B}$ be nearly subdirectly irreducible, with
$\alg{A}$ non-Boolean and
represented over $P$ with $d\deq \min(P\setminus M)$.
Then,
\begin{enumerate}
\item $Q(\alg{A})\subseteq Q(\alg{B})$ if and only if
there is a homomorphism $f\colon\alg{A}\to\alg{B}$ with
$f(\{d\})\neq 1$.
\item The inclusion problem
  $Q(\alg{A})\subseteq Q(\alg{B})$ is decidable.
\end{enumerate}  
\end{theorem}  

\begin{proof}
The backward direction of (1) is Lemma~\ref{lem:embedding}. For the forward direction,
if $\alg{A}\in Q(\alg{B})$, then by Lemma~\ref{lem:finite-products}
there is an embedding 
$e\colon \alg{A}\hookrightarrow \alg{B}^\ell$ for some
$\ell\in\mathbb{N}$. Therefore at least one projection $\pi$ must have
$(\pi\circ e)(\{d\})\neq 1$. Put $f \deq \pi\circ e$ to get the required
homomorphism.

For (2), if $\alg{A}$ is Boolean, then $\alg{A}\cong\alg 2^\ell$ and the inclusion
always holds. Otherwise it suffices to
check the finitely many maps from $A$ to $B$ for being homomorphisms not
sending $\{d\}$ to $1$.
\end{proof}

\begin{rem}
Our methods can also be applied to study the lattice of subvarieties
of $\var{Ho}$, dually isomorphic to \emph{axiomatic} extensions
of $\mathrm{IPC}_{\to,\neg}$, yielding analogous results. However
the lattice $\Lambda^v(\var{Ho})$ of \emph{subvarieties} of $\var{Ho}$
is much better known than $\Lambda^q(\var{Ho})$.
Several facts about it, such as its distributivity,
and the fact that completely join-irreducible members are generated
by finite subdirectly irreducible algebras, are folklore so we focussed on
$\Lambda^q(\var{Ho})$. One   
difference between $\Lambda^v(\var{Ho})$ and $\Lambda^q(\var{Ho})$
that we would like to point out here is that in $\Lambda^v(\var{Ho})$
any finitely generated variety has only finitely many subvarieties, whereas
it is not too difficult to find examples of quasivarieties
$Q(\alg{A})\in \Lambda^q(\var{Ho})$, generated by a single finite
algebra $\alg{A}$, even nearly subdirectly irreducible,
with an infinite chain of subquasivarieties. 
\end{rem}

\subsection{Logical reading}\label{subsec:logic1}
For a nearly subdirectly irreducible $\alg{A}$ let $\Log(\alg{A})$ denote
the extension of $\mathrm{IPC}_{\to,\neg}$ by all rules valid in
$\alg{A}$; under the dual isomorphism of Section~\ref{sec:intro},
$\Log(\alg{A})$ corresponds to $Q(\alg{A})$. Theorem~\ref{thm:structure}
then says that every finitary extension of $\mathrm{IPC}_{\to,\neg}$ is the
intersection of the logics $\Log(\alg{A})$ it is contained in, where
$\alg{A}$ ranges over finite nearly subdirectly irreducible algebras,
and that the lattice of all finitary extensions of $\mathrm{IPC}_{\to,\neg}$,
ordered by inclusion, is isomorphic to the lattice of upsets of the
countable poset $\mathbb{P}$ (equivalently, to the lattice of downsets of
its order dual). The logics $\Log(\alg{A})$ are exactly the completely
meet-prime elements of this lattice, and by Corollary~\ref{thm:inclusion-charact}(2)
the inclusion between two of them is decidable.

\subsection{Comparison with the generic case}
Hyndman and Nation~\cite{HN18} investigate the structure of
$\Lambda^q(\mathcal{K})$ for an arbitrary locally finite quasivariety, mostly
for finite signatures, so a brief comparison is in order.
It is folklore 
that for any locally finite quasivariety $\mathcal{K}$ of finite signature
every completely join-irreducible element of  
$\Lambda^q(\mathcal{K})$ is generated by a single finite
critical algebra. Moreover, every element of $\Lambda^q(\mathcal{K})$
is a join of completely join-irreducible elements
(cf.~\cite[Corollary~12 and Theorem~13]{HN18}).

Our Lemmas~\ref{lem:join-prime} and~\ref{lem:prime-singly-gen} are
analogues for completely join-prime elements and nearly subdirectly irreducible
algebras. In general, completely join-prime elements are completely
join-irreducible; the converse holds in infinitely join-distributive
(\emph{a fortiori} in completely distributive) lattices, but not in general,
thus no analogue of
our Theorem~\ref{thm:structure} holds in general. More precisely,
a generic restatement of Lemma~\ref{lem:join-prime} with `critical' replacing
`nearly subdirectly irreducible' fails. On the other hand, in our case
completely join-prime and completely join-irreducible elements coincide,
so an analogue of Lemma~\ref{lem:prime-singly-gen} with `critical' replacing
`nearly subdirectly irreducible' holds. Note however that
not every nearly subdirectly irreducible algebra is critical: for example,
$\alg{A}\times\alg{2}$ with a non-Boolean subdirectly irreducible $\alg{A}$ is
nearly subdirectly irreducible but not critical. For
$\alg{2}$ is a retract of $\alg{A}$ and 
$\alg{A}$ embeds into $\alg{A}\times\alg{2}$ via the map
$a\mapsto (a,g(a))$ where $g:\alg{A}\to\alg{2}$ is an onto homomorphism
(obtained by composing the double negation map of
Remark~\ref{rem:boolean} with a projection).
A complete characterisation of critical algebras in $\var{Ho}$
would be of interest, but we do not attempt it here.

\section{Continua of subquasivarieties}\label{sec:continua}

Here we will construct the continua of quasivarieties announced in
Section~\ref{sec:intro}. The following quasiequation will play an important
auxiliary role:
\begin{equation}\label{eq:W2}\tag{$W'_2$}
\begin{split}
 x\ra\neg y =1\ \doublewedge\
 (\neg\neg x\ra x)\ra z=1\ \doublewedge\
 (\neg\neg y\ra y)\ra z=1
 \quad\Ra\quad z=1,
\end{split}
\end{equation}
where $\neg x\deq x\ra0$. Read as an
inference rule, \eqref{eq:W2} says: infer $z$ from $x\ra\neg y$,
$(\neg\neg x\ra x)\ra z$ and $(\neg\neg y\ra y)\ra z$.
But we begin with some facts about heights that the
construction requires.

\subsection{Heights of posets and algebras}\label{subsec:heights}
For a poset $P$, define its height putting
\[
 h(P)=\sup\{|C|:C\text{ is a finite chain in }P\}.
\]
For an algebra $\alg A$, however, its height is 
$h(\Cm{\alg A})$. The discrepancy between the height of $\alg{A}$ in this
commonly used technical sense 
and the height of $\alg{A}$ as a poset comes from the
representation and is also well motivated equationally. 
For recall from~\cite{KSfree} the following terms:
\[
 d_1(x_1)=x_1,\qquad
 d_2(x_1,x_2)=((x_2\ra x_1)\ra x_2)\ra x_2,
\]
and $d_{k+1}=d_2(d_k,x_{k+1})$ for $k\geq2$. For $k\geq1$, let
\begin{equation}\label{eq:height}\tag{$H_k$}
 d_{k+1}(x_1,\ldots,x_{k+1})=1,
\end{equation}
and define
\[
 \var{Ho}_k\deq\{\alg A\in\var{Ho}:\alg A\models H_k\}.
\]
By \cite[Theorem~8.4]{KSfree}, for every $\alg A\in\var{Ho}$,
\[
 \alg A\in\var{Ho}_k
 \quad\Longleftrightarrow\quad h(\Cm{\alg A})\leq k.
\]
Thus, finite bounds on the heights of Hilbert algebras, and indeed  
of Heyting algebras, are equationally defined.   
In particular, $\var{Ho}_1$ is term-equivalent to the variety of Boolean algebras.
By Lemma~\ref{lem:height-poset}, if $\alg A$ is represented over $P$,
then $h(\Cm{\alg A})=h(P)$; this is how heights will be computed below.

\subsection{The algebras $\alg S^{(n)}$ and $\alg T^{(n)}$}

Fix $n\geq2$. Let $P^{(n)}$ be the poset on the set
$\{m,a,b,c,e,s_0,s_1,\ldots,s_{n-2}\}$
whose order is generated by taking $m$ as an isolated point, and
\[
 s_{n-2}<\cdots<s_1<s_0<b<a,\qquad s_0<e<c.
\]
So $M=\max P^{(n)}=\{m,a,c\}$,
and $h(P^{(n)}) = n+1$ (see
Figure~\ref{fig:witness} below). Define
\begin{equation}\label{eq:Sn}
 S^{(n)}=\pw(M)\cup\{\{p\}:p\in P^{(n)}\}
                     \cup\{\{a,e\},\{c,b\}\},
\end{equation}
and note that it is a hereditary family of antichains.  
By Lemma~\ref{lem:selection} the algebra
$\alg S^{(n)} = (S^{(n)};\ra,\emptyset,M)$ belongs to $\var{Ho}$,
and Lemma~\ref{lem:height-poset} gives
$h(\Cm{\alg S^{(n)}})=n+1$. It is clear by its construction that
$\alg{S}^{(n)}$ is nearly subdirectly irreducible.
Counting its elements gives
$|S^{(n)}| = n + 11$ (it will come in handy in the proof of Theorem~\ref{thm:main}). 

\begin{figure}[htbp]
\centering
\begin{tikzpicture}[x=.95cm,y=.85cm]
\foreach \shift/\name in {0/P,5.4/Q}{
 \begin{scope}[xshift=\shift cm]
 \node[bdot,label=above:{$a$}] (\name a) at (-.65,3.7) {};
 \node[bdot,label=above:{$c$}] (\name c) at (.65,3.7) {};
 \node[bdot,label=left:{$b$}] (\name b) at (-.65,2.65) {};
 \node[bdot,label=right:{$e$}] (\name d) at (.65,2.65) {};
 \node[bdot,label=right:{$s_0$}] (\name e) at (0,1.8) {};
 \draw[edge] (\name a)--(\name b)--(\name e)--(\name d)--(\name c);
 \draw[edge] (0,1.8)--(0,1.45);
 \draw[edge,densely dotted] (0,1.45)--(0,.35);
 \node[draw,fill=white,rounded corners,inner sep=4pt,align=center,font=\small]
    at (0,.55) {$s_1,\ldots,s_{n-2}$\\chain};
 \end{scope}
}
\node[bdot,label=above:{$m$}] at (-1.9,3.7) {};
\draw[-{Stealth[length=2mm]}] (1.6,2.4)--(3.9,2.4)
 node[midway,above,font=\small] {delete $m$};
\node at (-.4,-.25) {$P^{(n)}$};
\node at (5.4,-.25) {$P^{(n)}\setminus\{m\}$};
\end{tikzpicture}
\caption{The posets representing $\alg S^{(n)}$ and its quotient
$\alg T^{(n)}$. The box denotes a chain of $n-2$ points below $s_0$;
it is empty for $n=2$. The isolated point $m$ is deleted in the quotient.}
\label{fig:witness}
\end{figure}
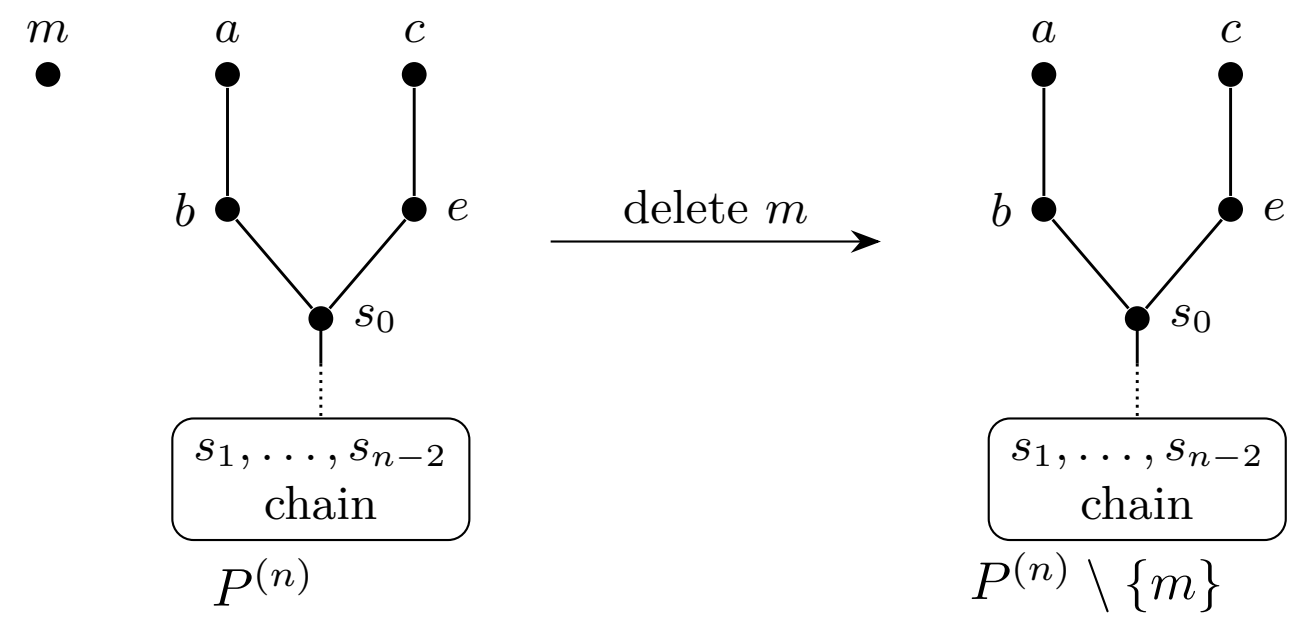

\begin{prop}\label{prop:Sn-W2}
For every $n\geq2$ we have $\alg S^{(n)}\models\eqref{eq:W2}$.
\end{prop}
\begin{proof}
Let $X,Y,R\in S^{(n)}$
satisfy the antecedent of~\eqref{eq:W2}. Its first conjunct gives
\[
 \emptyset=X\ra(Y\ra M)=M\setminus(X\cup Y),
\]
since the points of $M$ are maximal. Hence $m\in X$ or $m\in Y$ and by
symmetry we can assume $m\in X$. Every member of $S^{(n)}$ containing $m$ is a
subset of $M$, and every $U\subseteq M$ has $\neg\neg U=U$, so $\neg\neg X=X$.  
The second conjunct now yields
$\emptyset=(\neg\neg X\ra X)\ra R=\emptyset\ra R=R$.
Thus $R=\emptyset = 1$.
\end{proof}

Put $P'^{(n)}=P^{(n)}\setminus\{m\}$, $M' = M\setminus\{m\}$ and define
$T^{(n)}$ exactly as in~\eqref{eq:Sn}, with $P'^{(n)}$ in place of $P^{(n)}$
and $M'$ in place of $M$. Then let
\begin{equation}
 \alg T^{(n)}=(T^{(n)};\ra,\emptyset,\{a,c\}),
\end{equation}
where implication is computed in $P'^{(n)}$. By construction, $\alg{T}^{(n)}$ is
nearly subdirectly irreducible, indeed, subdirectly irreducible.

\begin{prop}\label{prop:quotient}
For every $n\geq2$, the algebra $\alg T^{(n)}$ is subdirectly irreducible
of height $n+1$ and is isomorphic to $\alg S^{(n)}/\theta_n$,
where $\theta_n=\theta^{\alg S^{(n)}}(1,\{m\})$. Moreover, 
$\alg T^{(n)}\not\models\eqref{eq:W2}$.
\end{prop}

\begin{proof}
The set $P'^{(n)}$ is an upset of $P^{(n)}$. By
Lemma~\ref{lem:restriction} (with $U=P'^{(n)}$),
\[
 r:S^{(n)}\longrightarrow T^{(n)},\qquad
 r(C)=C\cap P'^{(n)}=C\setminus\{m\}
\]
is a surjective homomorphism. Since $m$ is isolated,
$r(C)=\{m\}\ra C$, so~\eqref{eq:principal} gives
$\ker r = \theta^{\alg S^{(n)}}(1,\{m\}) =\theta_n$, so
$\alg S^{(n)}/\theta_n\cong\alg T^{(n)}$.
Lemma~\ref{lem:height-poset} gives $h(\Cm{\alg T^{(n)}})=n+1$.
The least point of $P'^{(n)}$ represents $\op$
in $\alg T^{(n)}$, so this algebra is subdirectly irreducible.

We verify the failure of~\eqref{eq:W2} directly in $\alg T^{(n)}$.
Put $X=\{a,e\}$, $Y=\{c,b\}$ and $R=\{s_0\}$. Here
$0=\{a,c\}$ and $1=\emptyset$, so
\[
 X\ra(Y\ra0)=X\ra\{a\}=\emptyset,
 \qquad \neg\neg X=\{a\},\quad \neg\neg Y=\{c\}.
\]
Furthermore,
\[
 (\neg\neg X\ra X)\ra R=\{e\}\ra\{s_0\}=\emptyset,
 \qquad
 (\neg\neg Y\ra Y)\ra R=\{b\}\ra\{s_0\}=\emptyset.
\]
All premises of~\eqref{eq:W2} hold, whereas $R=\{s_0\}\neq\emptyset=1$.
\end{proof}

For readers more familiar with Kripke frames it may be useful to remark that
the poset representation of $\alg{T}^{(n)}$ is a generated subframe of
the poset representation of $\alg{S}^{(n)}$. On the other hand,
$\alg{S}^{(n)}$ is not a direct product $\alg{T}^{(n)}\times \alg{2}$
(compare the cardinalities $n+11$ and $2(n+7)$),
nor is  $\alg{T}^{(n)}$ a subalgebra of $\alg{S}^{(n)}$, so the analogy
with Kripke frames is not perfect. The fact that
$\alg T^{(n)}$ falsifies $\eqref{eq:W2}$ can also be derived
from~\cite[Lemma~5.1]{CM10}.

\subsection{The algebras $\alg{E}_{n,k}$} 
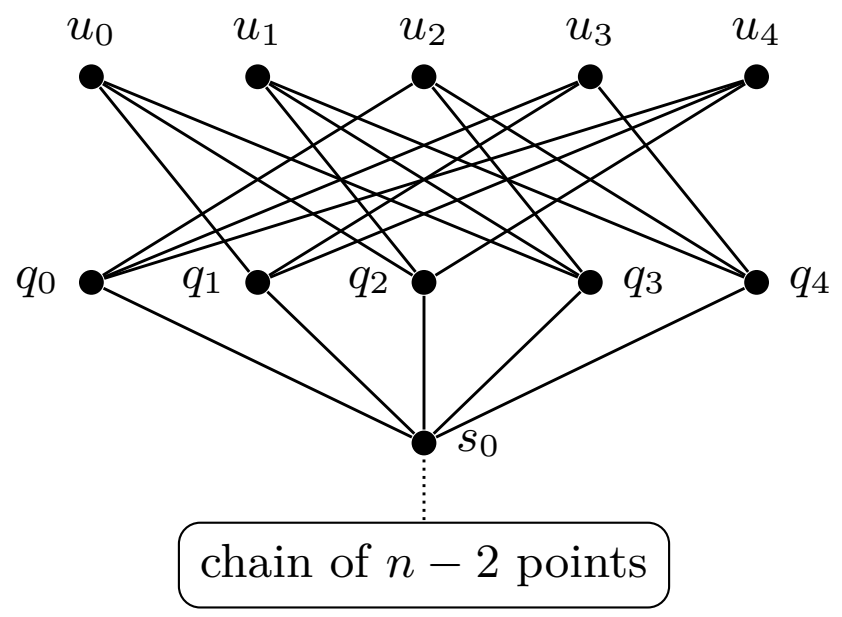
\begin{figure}[tbp]
\centering
\begin{tikzpicture}[x=1.1cm,y=.85cm]
\foreach \i in {0,1,2}{
 \node[bdot,label=left:{$q_{\i}$}] (q\i) at (\i,2.1) {};
}
\foreach \i in {3,4}{
 \node[bdot,label=right:{$q_{\i}$}] (q\i) at (\i,2.1) {};
}
\foreach \i in {0,1,2,3,4}{
  \node[bdot,label=above:{$u_{\i}$}] (u\i) at (\i,3.7) {};
}  
\foreach \i/\j in {0/2,0/3,0/4,1/0,1/3,1/4,2/0,2/1,2/4,3/0,3/1,3/2,4/1,4/2,4/3}{
 \draw[edge] (q\i)--(u\j);
}
\node[bdot,label=right:{$s_0$}] (r) at (2,.85) {};
\foreach \i in {0,1,2,3,4}{\draw[edge] (r)--(q\i);}
\draw[edge,densely dotted] (2,.85)--(2,.25);
\node[draw,fill=white,rounded corners,inner sep=4pt,align=center,font=\small]
 at (2,-.1) {chain of $n-2$ points};
\end{tikzpicture}
\caption{The poset $R_{n,5}$. The chain below $s_0$ 
is empty for $n=2$.}
\label{fig:family}
\end{figure}

For $n\geq2$ and $k\geq4$, let
$$
R_{n,k} = \{u_0,\ldots,u_{k-1},q_0,\ldots,q_{k-1},s_0,s_1,\ldots,s_{n-2}\}.
$$
Order $R_{n,k}$ as follows: the points $u_0,\ldots,u_{k-1}$ are maximal,
$s_{n-2}< \dots < s_1 < s_0$, $s_0 <q_i$ for all $i$, and 
$$
 q_i<u_j\quad\Longleftrightarrow\quad
 j\notin\{i,i+_k1\}
$$
where $+_k$ is addition modulo $k$ 
(see Figure~\ref{fig:family}).
Write $M_k=\{u_0,\ldots,u_{k-1}\}$ and define
$$
 E_{n,k}=\pw(M_k)\cup\{\{p\}:p\in R_{n,k}\},\qquad
 \alg E_{n,k}=(E_{n,k};\ra,\emptyset,M_k).
$$

\begin{lemma}\label{lem:E-properties}
For $n\geq2$ and $k\geq4$, the algebra $\alg E_{n,k}$, with $|E_{n,k}| =
2^k+k+n-1$, belongs to
$\var{Ho}$, has $h(\Cm{\alg E_{n,k}})=n+1$, is subdirectly irreducible
and satisfies~\eqref{eq:W2}. 
\end{lemma}

\begin{proof}
The family $E_{n,k}$ is nonempty, hereditary, and contains $M_k$, so
Lemma~\ref{lem:selection} gives $\alg E_{n,k}\in\var{Ho}$, and  simple counting
verifies its cardinality.

Since all singletons belong to $E_{n,k}$, by Lemma~\ref{lem:height-poset}
we get $h(\Cm{\alg E_{n,k}})=n+1$. Subdirect irreducibility is immediate by
noticing that $\op = \{s_{n-2}\}$.

It remains to show that $\alg E_{n,k}\models\eqref{eq:W2}$.
To this end, note that every element of $E_{n,k}$ is either a subset of $M_k$,
hence fixed by double negation, or a nonmaximal singleton whose negation is
$M_k$. If both $x$ and $y$ are 
of the latter kind, then $x\ra(y\ra M_k)= M_k\neq\emptyset = 1$.
Thus the first conjunct of~\eqref{eq:W2} forces one of $x,y$ to be
fixed by double negation. Without loss, assume
$\neg\neg x = x$. Then, the second conjunct forces $z=1$.
\end{proof}

\begin{lemma}\label{lem:nonembedding}
For fixed $n\geq2$ and $k,\ell\geq4$, the algebra $\alg E_{n,k}$
embeds into $\alg E_{n,\ell}$ if and only if $k=\ell$.
\end{lemma}

\begin{proof}
Suppose $f:\alg E_{n,k}\hookrightarrow\alg E_{n,\ell}$.
The elements fixed by double negation are precisely $\pw(M_k)$ and
$\pw(M_\ell)$, respectively. Hence $f$ restricts to an embedding of
$\pw(M_k)$ into $\pw(M_\ell)$, with $f(M_k)=M_\ell$. For $A,B\subseteq M_k$, we have
\[
 A\ra B=B\setminus A,\qquad (A\ra B)\ra B=A\cap B,
 \qquad\neg A=M_k\setminus A.
\]
It follows that $f$ preserves intersections and complements, therefore also
unions, on $\pw(M_k)$. Put
\[
 X_j=f(\{u_j\})\qquad(0\leq j<k).
\]
Since $f$ is injective and $f(\emptyset)=\emptyset$, every $X_j$ is
nonempty. For $i\neq j$,
\[
 X_i\cap X_j=f(\{u_i\}\cap\{u_j\})=\emptyset.
\]
By injectivity, the image of $\{q_i\}$ is not fixed by double negation;
hence it is a nonmaximal singleton. It cannot be $\{s_t\}$
for any $0\leq t\leq n-2$ since all these are comparable with every element of
the target algebra, whereas $\{q_i\}$ is incomparable with $\{q_j\}$ for $j\neq
i$; hence the image of $\{q_i\}$ must be some $\{q_j\}$. Since an embedding
preserves and reflects the order, there is an injective map $\sigma$ such that
\[
 f(\{q_i\})=\{q_{\sigma(i)}\}.
\]
Calculating implication for the two singletons of maximal points incomparable with $\{q_i\}$ gives
\[
 \{u_i\}\ra\{q_i\}=\{q_i\}
               =\{u_{i+_k1}\}\ra\{q_i\}.
\]
Applying $f$ yields
\[
X_i\ra \{q_{\sigma(i)}\}=\{q_{\sigma(i)}\}
               =X_{i+_k1}\ra\{q_{\sigma(i)}\},
\]
calculating implications we obtain
\[
\up{\{q_{\sigma(i)}\}}\cap X_i = \emptyset = \up{\{q_{\sigma(i)}\}}\cap X_{i+_k1}
\]
and therefore, by definition of order 
\[
 X_i\cup X_{i+_k1}\subseteq
       \{u_{\sigma(i)},u_{\sigma(i)+_\ell 1}\}.
\]
Since $X_i$ and $X_{i+_k1}$ are disjoint and nonempty, both are singletons, and
the same holds for every $i$. As $f$ preserves unions, we obtain
\[
 M_\ell=f(M_k)=f\left(\bigcup_{j=0}^{k-1}\{u_j\}\right)
              =\bigcup_{j=0}^{k-1}X_j.
\]
Since $|\bigcup_{j=0}^{k-1}X_j| = k$, we conclude that $k=\ell$.
The converse is immediate.
\end{proof}

Now we are ready to state and prove our second main result. 
Fix $n\geq2$. For $\Gamma\subseteq\Nstar$, define
\begin{equation}\label{eq:Qn}
 \mathcal Q^{(n)}_\Gamma=
 Q\bigl(\var{Ho}_n\cup\{\alg S^{(n)}\}
                   \cup\{\alg E_{n,k}:k\in \Gamma\}\bigr).
\end{equation}

\begin{theorem}\label{thm:main}
For every $n\geq2$, there is an order embedding
of $\pw(\mathbb N_{\geq4})$ into the interval $[\var{Ho}_n,\var{Ho}_{n+1}]$
in the subquasivariety lattice, given by
$\Gamma\longmapsto\mathcal Q^{(n)}_\Gamma$, such that
\[
\var{Ho}_n\subsetneq\mathcal Q^{(n)}_\Gamma \subsetneq\var{Ho}_{n+1}
\]
and no $\mathcal Q^{(n)}_\Gamma$ is a variety. Consequently, there are exactly
$2^{\aleph_0}$ quasivarieties which are not varieties and lie strictly
between $\var{Ho}_n$ and $\var{Ho}_{n+1}$.
\end{theorem}

\begin{proof}
Every generator of $\mathcal Q^{(n)}_\Gamma$ belongs to $\var{Ho}_{n+1}$, so
$\mathcal Q^{(n)}_\Gamma\subseteq\var{Ho}_{n+1}$. The inclusion
$\var{Ho}_n\subseteq\mathcal Q^{(n)}_\Gamma$ is strict, since
$\alg S^{(n)}\in\mathcal Q^{(n)}_\Gamma$ is of height $n+1$.

We claim that $\alg T^{(n)}\notin\mathcal Q^{(n)}_\Gamma$.
Otherwise, by Lemma~\ref{lem:finite-products} and
Corollary~\ref{cor:one-factor}, we would have
an embedding of $\alg T^{(n)}$ into one of the generators
of $\mathcal Q^{(n)}_\Gamma$. But $\alg T^{(n)}$ cannot embed into any member
of $\var{Ho}_n$, since $\alg T^{(n)}\not\models H_n$. 
Nor can it embed into $\alg S^{(n)}$ or any $\alg E_{n,k}$ as these
satisfy~\eqref{eq:W2}, and $\alg T^{(n)}$ does not.
This proves the claim. Since $\alg T^{(n)}\in\var{Ho}_{n+1}$, the
inclusion $\mathcal Q^{(n)}_\Gamma\subseteq\var{Ho}_{n+1}$ is strict.
Since $\alg T^{(n)}$ is a homomorphic image
of $\alg S^{(n)}$, no $\mathcal Q^{(n)}_\Gamma$ is a variety.

Next we show that
\begin{equation}\label{eq:recovery}
 \alg E_{n,k}\in\mathcal Q^{(n)}_\Gamma
             \quad\Longleftrightarrow\quad k\in \Gamma.
\end{equation}
The backward direction is immediate. For the forward direction,
finiteness and subdirect irreducibility, again by Corollary~\ref{cor:one-factor},
reduce membership to embedding into one generator. Height excludes a member of
$\var{Ho}_n$, while
\[
 |E_{n,k}|=2^k+k+n-1>n+11=|S^{(n)}|
\]
excludes $\alg S^{(n)}$; this inequality is the only place where the
restriction $k\geq4$ is used, as it fails for $k=3$.
Then Lemma~\ref{lem:nonembedding} leaves only a
generator with index $k$, proving $k\in \Gamma$.

Monotonicity of $Q$, together with~\eqref{eq:recovery}, gives
\[
 \mathcal Q^{(n)}_\Gamma\subseteq\mathcal Q^{(n)}_\Delta
                 \quad\Longleftrightarrow\quad \Gamma\subseteq \Delta,
\]
so there are at least $2^{\aleph_0}$ quasivarieties of the required kind.
There are at most $2^{\aleph_0}$ of them, since a quasivariety in a finite
signature is determined by the countable set of quasi-identities it
satisfies. 
\end{proof}

\begin{cor}\label{cor:continuum}
The variety $\var{Ho}$ has continuum many subquasivarieties which are
not varieties. Already $\var{Ho}_3$ contains continuum many such
quasivarieties, all strictly containing $\var{Ho}_2$.
\end{cor}
\begin{proof}
Take $n=2$ in Theorem~\ref{thm:main}.
\end{proof}

Theorem~\ref{thm:main} is false for $n=1$, since
$\var{Ho}_2$ is hereditarily structurally complete so
every subquasivariety of $\var{Ho}_2$ is a variety.
Hereditary structural completeness of $\var{Ho}_2$ was
noted in~\cite[Example~8.20]{KSfree} as an application
of~\cite[Theorem~8.19]{KSfree}, but it can also be derived
indirectly from~\cite[Theorem~5.3]{CM10}.

\subsection{Logical reading}\label{subsec:logic2}
For $k\geq1$ let $\mathrm{IPC}_{\to,\neg}+H_k$ denote the axiomatic
extension of $\mathrm{IPC}_{\to,\neg}$ by the formula
$d_{k+1}(x_1,\ldots,x_{k+1})$; it corresponds to the variety $\var{Ho}_k$,
and $\mathrm{IPC}_{\to,\neg}+H_1$ is classical logic in the language
$\{\to,\neg\}$. Theorem~\ref{thm:main} says that, for every $n\geq2$, there
are exactly $2^{\aleph_0}$ extensions of $\mathrm{IPC}_{\to,\neg}$ by
inference rules lying strictly between $\mathrm{IPC}_{\to,\neg}+H_{n+1}$
and $\mathrm{IPC}_{\to,\neg}+H_n$, none of which is axiomatic, and that
among them there is a copy of $\pw(\mathbb N_{\geq 4})$ ordered by
inclusion. Since there are only countably many finite sets of rules,
continuum many of these extensions are not finitely axiomatisable, even by
rules. By Corollary~\ref{cor:continuum}, all of this already happens
between $\mathrm{IPC}_{\to,\neg}+H_3$ and $\mathrm{IPC}_{\to,\neg}+H_2$,
whereas every extension of $\mathrm{IPC}_{\to,\neg}+H_2$ by rules is
axiomatic.

\subsection*{Acknowledgements}

Polish National Science Centre grant 2022/45/B/HS1/00606 is gratefully
acknowledged by both authors.

% Bibliography embedded to make this source self-contained.

\end{document}